\documentclass{amsart}
\usepackage[latin1]{inputenc}
\usepackage{amsmath}
\usepackage{amsthm}
\usepackage{amsfonts}
\usepackage{amssymb}
\newcommand{\ep}{\hspace*{\fill}$\Box$}

\begin{document}
\title[$P(f)=Q(g)$]{Rational decompositions of $p$-adic meromorphic functions}
\author[E. Mayerhofer]{Eberhard Mayerhofer}
\address{University of Vienna, Department of Mathematics, Nordbergstrasse 15, 1090 Vienna, Austria}
\email{eberhard.mayerhofer@univie.ac.at}
\begin{abstract} Let $K$ be a non
archimedean algebraically closed field of characteristic $\pi$,
complete for its ultrametric absolute value. In a recent paper by
Escassut and Yang (\cite{EY6}) polynomial decompositions
$P(f)=Q(g)$ for meromorphic functions $f$, $g$ on $K$ (resp. in a
disk $d(0, r^-)\subset K$) have been considered, and for a class
of polynomials $P$, $Q$, estimates for the Nevanlinna function
$T(\rho,f)$ have been derived.\\ In the present paper we consider
as a generalization rational decompositions of meromorphic
functions, i.e., we discuss properties of solutions $f$, $g$ of
the functional equation $P(f)=Q(g)$, where $P, \;Q$ are in $K(x)$
and satisfy a certain condition (M). We infer that in the case,
where $f$, $g$ are analytic functions, the Second Nevanlinna
Theorem yields an analogue result as in the mentioned paper
\cite{EY6}. However, if they are meromorphic, non trivial
estimates for $T(\rho,f)$ are more sophisticated.\end{abstract}
\maketitle
\section{Introduction} Throughout this paper, $K$ denotes an algebraically closed
ultrametric field of characteristic $\pi$, complete with respect
to the topology induced by its non-archi-\\medean valuation,
$K^*=K\setminus\{0\}$. Let ${\mathcal A}(K)$ denote the ring of
entire functions on $K$, and ${\mathcal M}(K)$ the field of
meromorphic functions in $K$, i.e., the field of fractions of
${\mathcal A}(K)$.\\ Moreover, for any real number $r>0$,
$d(a, r^-)=\{x:\vert x-a\vert<r\}$, i.e. the open ball with radius
$r>0$ and center $a\in K$; then similarly as above, ${\mathcal
A}(d(a, r^{-}))$ is the ring of analytic functions on $d(a, r^-)$,
i.e.: the ring of analytic functions with radius of convergence
$\rho\geq r$. The ring of meromorphic functions on $d(a, r^-)$ is
denoted by ${\mathcal M}(d(a, r^{-}))$.\\ We denote by
${\mathcal A}_b(d(a, r^-))$ the $K$-subalgebra of analytic
functions with bounded norm, furthermore ${\mathcal A}_u(d(a,
r^-))={\mathcal A}(d(a, r^-))\setminus{\mathcal A}_b(d(a, r^-))$.
Similarly, by \\ ${\mathcal M}_b(d(a, r^-))$ we denote the field of
fractions of ${\mathcal A}_b(d(a, r^-))$, and ${\mathcal M}_u(d(a,
r^-))={\mathcal M}(d(a, r^-))\setminus{\mathcal M}_b(d(a, r^-))$.
For $R>0$ we denote the interval $I$ as the set $I=[\rho, R[$,
where $0<\rho<R$, and for some $\rho>0$ we write
$J=[\rho,\infty[$.\medskip\\{\bf Notation in Nevanlinna Theory}
Let $R \in ]0,\infty[$, $f \in {\mathcal M}(d(0, r^{-}))$ (resp.
${\mathcal M}(K)$) such that $0$ is neither a zero nor a pole of
$f$. Let $w_{\alpha}(f)=n$ (resp. $w_{\alpha}(f)=-n$), if $f$ has
a zero (resp. a pole) of order n at $\alpha$. Then the functions
$Z$ and $N$ are defined as
\[
Z(\rho, f):=\sum_{w_{\alpha}>0,\vert \alpha\vert\leq \rho}
w_{\alpha}(f)\log \frac{\rho}{\vert \alpha \vert}
\]
and $N(\rho,f):=Z(\rho, 1/f)$, moreover the Nevanlinna function is
given by
\[
T(\rho,f):=\max\{Z(\rho, f), N(\rho,f)\}
\]
In addition we use similar functions not respecting multiplicities
of zeros (resp. poles):
\[
\widetilde Z(\rho, f):=\sum_{w_{\alpha}>0,\vert \alpha\vert\leq
\rho} \log \frac{\rho}{\vert \alpha \vert}
\]
and similarly $\widetilde N(\rho,f):=\widetilde Z(\rho, 1/f)$.
\medskip\\{\bf Notation in positive characteristic.} If $\pi\neq
0$ we define the characteristic exponent $\chi:=\pi$, otherwise we
set $\chi:=1$. Due to \cite{EY6} we call the ramification index of
$h$ in ${\mathcal M}(d(0, R^-))$(resp. ${\mathcal M}(K)$) the
unique integer $t$ such that $\sqrt[\chi^t]{h}$ belongs to
${\mathcal M}(d(0, R^-))$(resp. ${\mathcal M}(K)$). If $\pi=0$,
every function $h$ in ${\mathcal M}(d(0, R^-))$(resp. ${\mathcal
M}(K)$) has ramification index equals 0.\\ We note that in nonzero
characteristic, the counting functions $\widetilde N$ (resp. $
\widetilde Z$) of poles (resp. zeros) are defined slightly differently.
For more information, we refer to \cite{BE4}.
\\
In the present paper we apply the Second Nevanlinna Theorem due to
Boutabaa and Escassut (\cite{BE4}, Theorem 2):
\medskip\\{\bf Theorem N.} {\it Let $\alpha_1,\dots,\alpha_n\in
K$, with $n\geq 2$, and let $f\in{\mathcal M}(d(0, R^-))$ (resp.
$f\in {\mathcal M}(K)$) of ramification index $s$, have no zero
and no pole at $0$. Let $S:=\{
\sqrt[\chi^s]{\alpha_1},\dots,\sqrt[\chi^s]{\alpha_s} \}$. Assume
that f, $\sqrt[\chi^s]{f}$, $f-\alpha_j$ have no zero and no pole
at $0$ ($1\leq j\leq n$). Then we have:
\[
\frac{(n-1)T(r,f)}{\chi^s}\leq \sum_{i=1}^n\widetilde
Z(r,f-\alpha_i)+\widetilde N(r,f)-\log r+O(1),\;\; r\;\in\;
I\;\;( r\;\in\; J)\;
\]}
\\\\ Many Applications of the Nevanlinna Theory to Functional
and Differential Equations have been worked out in the last years,
and the Theory has only recently been generalized to fields of
characteristic $\pi$ (see \cite{BE4} which contains the Theorem from
above, resp. \cite{HY8}). One of the most famous examples, where
the archimedean Theorem is due to F. Gross \cite{G7} (generalizations were firstly made by N. Toda \cite{T9}, and the
non-archimedean work is due to A. Boutabaa \cite{B1}) is the
equation
\[
f^n+g^m=1
\]
A recent $p$-adic article on this topic deals with unbounded
meromorphic solutions in a ball (\cite{BE4}). In
characteristic zero, the most comprehensive work on this class of
functional equations can be found in \cite{BE2}. \\ Here we
discuss properties of analytic or meromorphic solutions $f$, $g$ of the functional equation
\begin{equation}\label{eq1}
P(f)=Q(g)
\end{equation}
where $P$ and $Q$ are certain rational functions on $K$.
\\\\
We are starting with analytic functions in section $2$ and receive similar conclusions as in a recent paper due to Escassut
and Yang (\cite{EY6}), where $P$, $Q$ are elements in $K[x]$. However,
the meromorphic case turns out to be more sophisticated than the
analytic one, i.e., it is more complicated to derive non trivial
estimations for the Nevanlinna function $T(\rho, f)$. This case
is worked out in section $3$.\\ \\When not
explicitly stated, we write $P=R/S$, $Q=V/W$, where $R$, $S$ ,$V$,
$W$ $\in K[x],\; (R,S):=\gcd(R,S)=1,\; (V,W)=1$, and the degrees
of $P$ (resp. $Q$) are defined by $p:=\deg P=\max\{\deg R, \deg
S\}$ (resp. $q:=\deg Q=\max\{\deg V, \deg W\}$).\medskip\\
 {\bf Preparatory statements}\medskip\\ {\bf Theorem 0.1.} {\it Let
$f\in\mathcal{A}(K)$ such that $f(K)\subset K^*$ (i.e. f has no
zero in K). Then f is a constant. }
\begin{proof}
Let $f(x)=\sum_{j\geq 0}a_jx^j$. It is well known that the number of
zeros of $f$ in any disk $d(0,r):=\{x\in K:\vert x\vert\leq r\}$
is equal to the largest integer k such that $\vert a_k\vert
r^k=\sup_{i\geq 0}\vert a_i\vert r^i$ (see Theorem 23.5 in
\cite{E5}). Hence, if $f$ has no
zero in $K$, obviously for all $n>0:\;a_n=0$.\end{proof}
{\bf Corollary 0.2.} {\it If $f$, $g$ $\in {\mathcal A}(K)$ have
the same zeros respecting multiplicity, then $\frac{f}{g}$ is a
constant.}\ep\medskip\\
{\bf Theorem 0.3.} (\cite{BE3}) {\it Let f be in ${\mathcal
M}(d(0, r^{-}))$ with $f(0)\neq 0,\infty$. Then $f$ belongs to
${\mathcal M}_b(d(0, r^{-}))$ if and only if $T(\rho,f)$ is
bounded in $[0,r[$.}\medskip\\
{\bf Theorem 0.4.} (\cite{BE3}) {\it Let $f$ b in ${\mathcal
M}(d(0, r^{-}))$ and let $P\in K[x]$ be of degree $n$. Then
$T(\rho,P(f))=nT(\rho,f)+O(1)$.}\medskip\\
{\bf Remark 0.5.} (\cite{BE4}, \cite{EY6}) There is a well defined
mapping
\[
\sqrt[\chi]{}:\quad\quad x\mapsto\sqrt[\chi]{x}
\]
which is a homomorphism on $K$ and can be extended to a
homomorphism on $K[x]$ (or even to one on $K(x)$) in the following
way: Let $R=\sum_{i=0}^n a_i
x^i=\mu(x-{\alpha_1})\dots(x-{\alpha_n})$ in $K[x]$, then we
define
\[
\sqrt[\chi]{}:\quad\quad \sqrt[\chi]{P}(x):=\sum_{i=0}^n
\sqrt[\chi]{a_i}x^i=\sqrt[\chi]{\mu}(x-\sqrt[\chi]{\alpha_1})\dots(x-\sqrt[\chi]{\alpha_n})
\]
{\bf Lemma 0.6.} (\cite{BE2}, \cite{BE4}) {\it Suppose $\pi\neq
0$. Let $r>0$ and $f\in{\mathcal M}(d(0,r^-))$(resp.
$f\in{\mathcal M}(K)$). Then $\sqrt[\chi]{f}$ belongs to
${\mathcal M}(d(0,r^-))$(resp. ${\mathcal M}(K)$) if and only if
$f'\equiv 0$. Moreover, there exists a unique $t\in\mathbb{N}$
such that $\sqrt[\chi^t]{f}$ in ${\mathcal M}(d(0,r^-))$(resp.
${\mathcal M}(K)$) and $(\sqrt[\chi^t]{f})'\neq 0$.}
\section{Decompositions of analytic functions}
First, we apply Theorem 2.9 (\cite{EY6}) to the question
(\ref{eq1}) for entire $f$, $g$. In section 2.2
we prove a somewhat similar result with appropriate
conditions, particularly tailored to our "rational" problem.
Indeed this does not only yield a quite more general result (see
examples 2.2.6, 2.2.8), but also an analogue result for elements
$f$, $g$ in ${\mathcal A}(d(0,r^-))$.
\subsection{An Application of a previous paper}
(\cite{EY6})\medskip\\
{\bf Remark 2.1.1.} Let $f$, $g$ be in $A(K)$ (resp $A(d(0,
r^-))$), satisfying (\ref{eq1}). Any pole $b$ of $P(f)$ is a zero
of $S(f)$, hence a zero of $W(g)$ of the same order. In the
meromorphic case, this conclusion is wrong.\medskip\\
{\bf Theorem 2.1.2.} {\it Let $f$, $g$ $\in{\mathcal A}(K)$ solve
(\ref{eq1}), then there exists a constant $\lambda\in K^*$ such
that
\[
S(f)=\lambda W(g)
\]}\begin{proof}
Since $f$, $g$ are analytic functions and $S(f)$ and
$W(g)$ have the same zeros of the same order, $S(f)/W(g)$ is a
constant by Corollary 0.2.\end{proof}
Thus our problems reads
\begin{equation}\label{eq25} R(f)=\lambda
S(g),\;\lambda\in K^*
\end{equation}
Although $\lambda$ is an undetermined constant, we are able to
apply the following Theorem (Theorem 2.9 in \cite{EY6}) to it,
worked out for decompositions $A(f)=B(g)$ for $f$, $g$, where $A$,
$B$ are polynomials (in order to avoid confusion, we write $A$,
$B$ instead of $P$, $Q$ used in \cite{EY6}):\medskip\\ {\bf
Theorem 2.1.3.} (\cite{EY6}) {\it Let $A$, $B$ be in $K[x]$ with
$A'B'$ not identically zero, such that \\ $2\leq\min\{\deg
A, \deg B\}$. Assume that there exist $k$ distinct zeros
$c_1,\dots,c_k$ of $A'$ such that $A(c_i)\neq A(c_j)\forall i\neq
j$ and $A(c_i)\neq B(d)$ for every zero d of $B'$($i=1,\dots,k$).
Assume that there exist two nonconstant functions $f$, $g$
$\in{\mathcal M}(K)$ such that $A(f)=B(g)$, and let $t=\nu(f)$.
Then $q\leq p$ and f satisfies
\[
\widetilde N(\rho,f)\geq \frac{k\deg B-\deg A}{\chi^t \deg
B}T(\rho,f)+\log \rho+O(1)
\]
Moreover, if $\frac{p}{2}<q$, then $k=1$ and $c_1$ is a simple
root of $A'$.}\medskip\\
{\bf Theorem 2.1.4.} {\it Let $P=R/S$, $Q=V/W$, where $R$, $S$,
$V$, $W$ are in $K[x]$, $(R,S)=1$, $(V,W)=1$. Let $c_1,\dots,c_k$
be zeros of $R'$ such that $R(c_i)\neq R(c_j)\forall i\neq j$.
Moreover, let the degree of $V$ satisfy $l=k-\deg V+1>0$. Then, if
$f,g \in{\mathcal A}(K)\setminus K$ solve $P(f)=Q(g)$, we have
\[
0\geq \frac{l\deg R-\deg V}{\chi^t\deg V}T(\rho, f)+\log \rho+O(1)
\]
i.e.: $\deg V(\deg R+1)>(k+1)\deg R$. }\begin{proof}
Due to Theorem 2.1.2 we have $R(f)=\lambda V(g)$ for some $\lambda\in
K^*$; set $A:=R,\; B:=\lambda V$. Now, obviously there exist at
least $l$ distinct roots $c_{j_1},\dots,c_{j_l}$ of $A'$, such
that $A(c_{j_r})\neq B(d)$ for any zero $d$ of $B$ ($r=1,\dots,l$
with $l=k-\deg V+1>0$), and $\forall i\neq j, (i,j \in
\{j_1,\dots,j_l\})\;:A(c_i)\neq A(c_j)$. Trivially, $\widetilde
N(\rho,f)$ is identically zero for non constant analytic
$f$.\end{proof}
{\bf Remark 2.1.5.} Unfortunately the condition
of Theorem 2.1.4 implies \\ $\deg V\leq\deg R-1$, which
follows from $\deg V\leq k<\deg R$, since $k\leq \deg R'$. This
inequality together with the statement of Theorem 2.1.4 tells us
therefore:
\[
\frac{(k+1)\deg R}{\deg R+1}<\deg V<\deg R
\]
In the next section, however, we present conditions for $P$, $Q$,
where not necessarily  $\deg V<\deg R$, such that (\ref{eq1}) has
only non constant entire solutions. Also, we consider the case
where $f$, $g$ are unbounded analytic functions inside a disk
$d(0,r^-)$. \medskip\\{\bf Corollary 2.1.6.}  {\it Let $P,\ Q$ be
in $ K(x)$ with $P'Q'$ not identically zero and let $p=\deg(R),\
q=\deg(V)$ with $2\leq \min(p,q)$ and $\frac{p}{2}<q$. Assume that
there exist $q$ distinct zeros $c_i,\; (1\leq i\leq q)$ of $R'$
such that $R(c_i)\neq R(c_j)\; \forall i\neq j$. If two functions
$f,\ g\in {\mathcal A}(K)$ satisfy $P(f)=Q(g)$, then $f$ and $g$
are constants. } \begin{proof}
Assume that two
functions $f,\ g\in {\mathcal A}(K)$ satisfy $P(f)=Q(g)$. By
Theorem 2.1.2 there exists $\lambda \in K$ such that $R(f)=\lambda
V(g)$. Let $d_1,\dots,d_n$ be the distinct zeros of $V'$. We
notice that $n\leq q-1$. In order to apply Theorem 2.10 in
\cite{EY6}, we only have to check that there exists a zero $c_k$
of $R'$ satisfying $R(c_k)\neq \lambda V(d_j)$ for every
$j=1,\dots,n$. Suppose it is not true. Then, up to a reordering,
we can assume that $R(c_1)=\lambda V(d_1),\dots, R(c_1)=\lambda
V(d_n)$. Since $q>n$ and since $R(c_i)\neq R(c_j)\; \forall i\neq
j$, we then have $R(c_q)\neq \lambda V(d_j)\; \forall
j=1,\dots,n$. Thus, we can apply Corollary 2.10 in \cite{EY6} to
the polynomial $A:=R$ and $B:=\lambda V$.\end{proof}
\subsection{Generalizations of} \cite{EY6}\medskip\\
{\bf Remark 2.2.1.} We further distinguish some cases:
Suppose $f$, $g$ are non constant entire solutions (resp.
unbounded analytic solutions in a disk $d(a,r^-)$) of (\ref{eq1}),
then by using growth at infinity (resp. growth, when
$\rho\rightarrow r-$) we have
\\{\bf Case 1.} $\deg V>\deg W$, then obviously  $\deg R>\deg S$,
\\{\bf Case 2.} $\deg V<\deg W$, then obviously  $\deg R<\deg S$,
\\{\bf Case 3.} $\deg V=\deg W$, then obviously  $\deg R=\deg S$,\\
furthermore we obtain in any case, when $\rho\rightarrow \infty$
(resp. when $\rho\rightarrow r-$ )
\begin{equation}\label{eq2}
p T(\rho,f)=q T(\rho, g)+ O(1)
\end{equation}
which follows from the functional equation (1): \\ $T(\rho,
P(f))=\max\{Z(\rho,P(f)),\; N(\rho,P(f))\}$, and since
$Z(\rho,P(f))=Z(\rho, R(f))$, \\$N(\rho,P(f))=Z(\rho, S(f))$, we
have by Theorem 0.4, $T(\rho, P(f))=\max\{r,s\}T(\rho,f)+O(1)=p
T(\rho,f)+O(1)$.\medskip\\ In the present paper our statements on
rational decompositions of meromorphic functions always concern a
specific class of rational functions $P$, $ Q$, admitting certain
decompositions themselves (see Lemma 2.2.2, below). They are described by the following condition to which we always refer:\medskip\\
{\bf Condition (M)} Let $P, Q\in K(x)$ and denote the zeros of
$P'$ by $c_1,\dots c_k$.\\ $P$, $Q$ are said to satisfy
Condition (M), if
\begin{enumerate}
\item $P'Q'\neq 0$

\item $P=R/S$, $(R,S)=1$, $Q=V/W$, $V$, $W$ monic, $(V,W)=1$ ($R$,
$S$, $V$, $W$ $\in K[x]$)\item $k>0$, and for any $i$
$\in\{1,\dots, k\}$ we have
\[
\left (Q(d)\neq {P(c_i)}\right )\quad \wedge\quad\left( W(d)\neq
0\right),
\]
for any zero $d$ of $V'-W'P(c_i)$, \item $P(c_i)\neq P(c_j)$
whenever $i \neq j$, \item Finally, if $v=w$ we assume $\forall
i\in\{1,\dots,k\}: P(c_i)\neq 1$.
\end{enumerate}{\bf Remark on Condition (M).}
Let $f$, $g$ be non constant entire functions. If we set $S\in
K^*$, it easily follows that also $W\in K^*$, moreover $P, Q\in
K[x]$ and (like before, denote the zeros of $P'$ by $c_1,\dots
c_k$) satisfy:
\begin{enumerate}
\item $P'Q'\neq 0$
\item $k>0$, and for any $i$ $\in\{1,\dots, k\}$ we have $
Q(d)\neq {P(c_i)}$ for any zero $d$ of $Q'$, \item $P(c_i)\neq
P(c_j)$ whenever $i \neq j$
\end{enumerate}
On this condition Theorem 2.1 and 2.9 are based in \cite{EY6}. In
this sense, our paper can be considered as a generalization to
part of \cite{EY6}. \\ Now -we are ready to state the basic Lemma:
\medskip\\ {\bf Lemma 2.2.2.} {\it Let $P$, $Q$ satisfy Condition
(M). Then for any $i\in{1,\dots,k}$, $P-P(c_i)$ resp.
$(Q-P(c_i))W$ we have the following factorizations:
\begin{equation}\label{eq3}
P(x)-P(c_i)=(x-c_i)^{s_i}R_i(x),\quad s_i\geq 2,\quad R_i(c_i)\neq
0
\end{equation}
and
\begin{equation}\label{eq4}
(Q(x)-P(c_i))W(x)=\prod_{j=1}^q (x-b_{i,j})
\end{equation}
Furthermore, the set $\{b_{i,j}\}$ consists of $qk$ distinct
elements.}\begin{proof} 
Let $i\in\{1,\dots,k\}$ be arbitrary, but fixed. Since $P'(c_i)=0$ we can clearly write
$P(x)-P(c_i)=(x-c_i)^{s_i}R_i(x)$, with $R_i(c_i)\neq 0$, and
$s_i\geq 2$; indeed, suppose $s_i=1$, then for the derivative we
have $P'(x)=(x-c_i)R_i'(x)+R_i(x)$, so that $P'(c_i)=R_i(c_i)$,
which is a contradiction. \\ In any case, we have
$\deg((Q-P(c_i))W)=\deg(Q)=\max\{\deg V, \deg W\}=q$: In Case 1
and Case 2 this is obvious, and in Case 3 we infer this by the
additional condition $P(c_i)\neq 1$ and $V$, $W$ being monic
polynomials. Thus we can write $(Q(x)-P(c_i))W(x)=\prod_{j=1}^q
(x-b_{i,j})$; furthermore for any fixed $i$, $b_{i,j}\neq b_{i,j'}$,
since for any $d$ with $W(d)=0$ or $Q(d)-P(c_i)=0$ we have
$V'(d)-W'(d)P(c_i)\neq 0$.\\ Now, let $(i,j)\neq (i',j')$, then
$\varphi(x):=(Q(x)-P(c_i))-(Q(x)-P(c_{i'}))=P(c_i)-P(c_{i'})\neq
0$ is a constant function different from zero; on the other hand,
assume $b_{i,j}= b_{i',j'}$, then by the right side of
decomposition (\ref{eq4}) we infer $\varphi(x)=0$, which
contradicts our assumption and thus $b_{i,j}\neq
b_{i',j'}$.\end{proof} Let $R_i$ ($i=1,\dots,k$) be the rational
functions due to the notation of Lemma 2.2.2. The following lemma
presents upper bounds for their degrees ($i=1,\dots,k$):\medskip\\
{\bf Lemma 2.2.3.} {\it Let $P$, $Q$ satisfy condition (M). Then
for any $i\in\{1,\dots,k\}$, $R_i$ can be written in the way
$R_i=\frac{A_i}{B_i}$ with $(A_i,B_i)=1$, where $B_i=S$. Moreover
we have $\deg A_i\leq\max\{\deg R,\deg S\}-s_i=p-s_i$.}
\begin{proof} We may write $R_i=\frac{A_i}{B_i}$ and from
(\ref{eq4}) we get
\[
    (x-c_i)^{s_i}\frac{A_i(x)}{B_i(x)}=\frac{R(x)-P(c_i)S(x)}{S(x)}
\]
$(R,S)=1$ clearly implies $(R(x)-P(c_i)S(x), S(x))=1$, thus
$B_i=S$ and $\deg A_i=\deg(R(x)-P(c_i)S(x))-s_i$.
\\{\bf Case 1.} $\deg(R(x)-P(c_i)S(x))=r$, i.e.
    $\deg A_i=\deg(R(x)-P(c_i)S(x))-s_i=r-s_i=p-s_i$.
\\{\bf Case 2.} We get exactly in the same way as before,
    switching the roles of $r$ and $s$: $\deg(R(x)-P(c_i)S(x))=s$, i.e.
    $\deg A_i=\deg(R(x)-P(c_i)S(x))-s_i=s-s_i=p-s_i$.
\\{\bf Case 3.} Obviously $R$ might be not monic; we
    conclude $\deg A_i=\deg(R(x)-P(c_i)S(x))-s_i\leq \deg S-s_i=\deg R-s_i=\deg P-s_i$.
\end{proof}
Particularly for the case where $K$ has nonzero characteristic, we note two useful Lemmas:\medskip\\{\bf Lemma $\Pi_1$.} {\it If for
$P$, $Q$ in $K(x)$, $P'Q'\neq 0$ and $f$, $g$ in $\in{\mathcal
M}(K)\setminus K$ (resp. ${\mathcal M}_u(d(0,r^-))$) satisfy
(\ref{eq1}), then $f'\equiv 0\Leftrightarrow g'\equiv 0$.
Obviously, $t=\nu(f)=\nu(g)$.} \begin{proof} Say $f'=0$,
then by the derivative of (\ref{eq1}) we can see that $g'=0$:
Firstly, $\mathcal M(K)$ (resp. $\mathcal M(d(0,r^{-}))$) is a
field, secondly $Q'$ is not identically zero by our assumption,
i.e. it vanishes at finitely many points only; and since $g$ takes
infinitely many values, $Q'(g)$ is not identically zero.
Conversely, by the same argument if $f'\neq 0$, then $g'\neq
0$.\end{proof}
{\bf Lemma $\Pi_2$.} {\it Let $P$, $Q$ satisfy Condition (M) and
$f'\equiv 0$. Then Condition (M) is satisfied by
$P_1:=\sqrt[\chi]{P}$, $Q_1:=\sqrt[\chi]{Q}$. Moreover if $f,\;g\;
\in\;\mathcal{M}(K)$ (resp. ${\mathcal M}(d(0,r^-))$) satisfy
(\ref{eq1}), then $P_1(f_1)=Q_1(g_1)$, where
$f_1:=\sqrt[\chi]{f}$, $g_1:=\sqrt[\chi]{g}$. Thus by repeating
the same process $t$ times, where $t$ is the unique integer from
Lemma 0.6, we derive $P_t(f_t)=Q_t(g_t)$ and Condition (M) is
satisfied by. $P_t, Q_t$ (where we denote similarly
$f_t:=\sqrt[\chi^t]{f}$, $g_t:=\sqrt[\chi^t]{g}$).}\medskip\\
\begin{proof} Use Remark 0.5. \end{proof}
{\bf Theorem 2.2.4.} {\it Let $f$, $g$ $\in{\mathcal A}(K)$ be non
constant solutions of (\ref{eq1}), where $P$, $Q$ satisfy
condition (M) and let $t:=\nu(f)=\nu(g)$. Then
\begin{equation}\label{eq40}
0\geq (\frac{kq-p}{\chi^tq}) T(\rho, f)+\log \rho+O(1),
\end{equation} i.e., $qk-p<0$.}\medskip\\
{\bf Theorem 2.2.5.} {\it Let $f$, $g$ $\in{\mathcal
A}_u(d(0,r^-))$ be solutions to (\ref{eq1}), where $P$, $Q$
satisfy condition (M) and let $t:=\nu(f)=\nu(g)$. Then
\begin{equation}\label{eq24}
0\geq (\frac{kq-p}{\chi^tq}) T(\rho, f)+O(1),
\end{equation}
i.e., $qk-p \leq 0$.} \medskip\\ For the proof of these two
theorems we refer to section 4, where a unified proof including
analogue statements on meromorphic functions (see section 3) is
given. \subsection{Examples} Let $\deg R=\deg S=\deg V=\deg W=2$,
then due to Remark 2.1.5 Theorem 2.1.4 can not be applied. But
2.2.4 works:  To demonstrate this explicitly, we consider the
following example:\medskip\\{\bf Example 2.2.6.} Let $K=\mathbb
C_p$, let $P(x)=\frac{R(x)}{x^2}$, $\deg R=2$,
$Q(x)=\frac{V(x)}{W(x)}=\frac{x^2}{x^2+x+1}$. Write
$R(x)=ax^2+bx+c$, $a\neq 0$, $b\neq 0$, $\frac{a-c}{b}=\sqrt 3$.
Then
\[
P'(x)=\frac{-bx^2+x(2a-2c)+b}{(x^2+1)^2}=\frac{-bx^2+2(b\sqrt
3)x+b}{(x^2+1)^2}=-b\frac{x^2-2\sqrt 3x-1}{(x^2+1)^2}
\]
and $P'$ has two distinct roots
\[
c_1=\sqrt 3+2, \quad c_2=\sqrt 3-2
\]
which yields
\[
P(c_1)=\frac{(4\sqrt 3+7)b+(4+2\sqrt 3)c}{4+2\sqrt 3},\quad
P(c_2)=\frac{(4\sqrt 3-7)b+(4-2\sqrt 3)c}{4-2\sqrt 3}
\]
We may set $P(c_1)=\frac{2}{3}$, $P(c_2)=-\frac{4}{3}$. Thus, the
only zero $d$ of $V'-P(c_1)W'=x(2-2P(c_1))-P(c_1)$ is
$d=\frac{P(c_1)}{2(1-P(c_1))}=1$, furthermore $W(d)=W(1)\neq 0$
and $1 \neq \frac{2}{3}=P(c_1)\neq Q(d)=\frac{1}{3}$; similarly,
the only zero $d$ of $V'-P(c_2)W'=x(2-2P(c_2))-P(c_2)$ is
$d=\frac{P(c_2)}{2(1-P(c_2))}=-\frac{2}{7}$, furthermore
$W(d)=W(-\frac{2}{7})\neq 0$ and $1\neq -\frac{4}{3}=P(c_2)\neq
Q(d)=\frac{4}{39}$.\\ By applying Theorem 2.2.4, we conclude that
there are no non constant, entire solutions $f$, $g$ of the
equation
\[
\frac{(\frac{\sqrt 3}{2}-\frac{1}{3})f^2+f-(\frac{\sqrt
3}{2}+\frac{1}{3})}{f^2}=\frac{g^2}{g^2+g+1}
\]
Note that in the case $b=0$ the problem turns out to be almost
trivial (use Theorem 2.1.2).\\ Theorem 2.2.5 can be applied, too, since
$kq-p=4-2=2$, and we get, that there are no
unbounded elements $f$, $g$ in ${\mathcal A}(d(0,r^-))$ having the
above decomposition.\medskip\\ Example 2.2.5 shows that in the
following case condition (M) does not yield an empty set of
rational functions $P$, $Q$:\medskip\\{\bf Corollary 2.2.7.} {\it
Let $P=R/S$, $Q=V/W$, where $R$, $S$, $V$, $W$ are polynomials
over K, V, W, monic, (R,S)$=1$, (V,W)$=1$, each of which having
degree two. Let $P'$ have two distinct zeros $c_i\;(i=1,2)$ such
that $P(c_1)\neq P(c_2)$ and $P(c_i)\neq 1\;(i=1,2)$. Assume that
for any i, $\deg (V'-P(c_i)W')>0$ and let $d_i$ be its unique
zero. Also, suppose $P(c_i)\neq Q(d_i)$ and $W(d_i)\neq
0\;(i=1,2)$. Then for a pair $(f,g)\in {\mathcal
A}(K)\times{\mathcal A}(K)$ having the decomposition $P(f)=Q(g)$
it follows that $(f,g)\in K^2$.}\ep \medskip\\ Similarly to
Example 2.2.6 we show now, that there exist rational functions
$P$, $Q$, $\deg R=\deg V=2, \deg S=2$ which satisfy Condition
(M).\medskip\\{\bf Example 2.2.8.} Let $K=\mathbb C_p$. For $f$,
$g$ $\in {\mathcal A}(K) $(resp. ${\mathcal A}(d(0,r^-))$), we
consider the functional equation
\[
\frac{a f^2+b f+c}{f^3}=\frac{g^2}{g^3-6g^2+11 g+6}
\]
where $a$, $b$, $c$ in $K$ are chosen in a way, that Condition (M)
is satisfied: We write $R=ax^2+bx+c$, $S(x)=x^3$, $V(x)=x^2$,
$W(x)=(x-1)(x-2)(x-3)$. Clearly we have to choose $c\neq 0$, such
that $(R(x),x^3)=1$. Whenever $a\neq 0$, the derivative $P'$ is
\[
P'(x)=\frac{(-a(x^2+2b/a x+3c/a))}{x^4}
\]
For $c=\frac{b^2}{3a}$, $P'$ has a single zero of multiplicity two
only: $c_1=-\frac{b}{a}$; furthermore $P(c_1)=-\frac{a^2}{3b}$. We
set $t:=-P(c_1)=\frac{a^2}{3b}$. Now the reader can easily verify
that we may choose $t$ in such a way that
\begin{enumerate}
\item $V'-P(c_1)W'=3tx^2+(2-12t)x+11t=0$ has one single solution
of multiplicity two, $d=\frac{6t-1}{3t}$, \item $d\notin
\{1,2,3\}$ (i.e., $W(d)\neq 0$) \item $P(c_1)\neq Q(d)$.
\end{enumerate}
Since $p=q=3,\; k=1$, Theorem 2.2.4 assures us that there are no
non constant entire functions $f$, $g$ satisfying the functional
equation from above. However, elements in ${\mathcal
A}_u(d(0,r^-))$ with this specific decomposition might
exist.\medskip\\{\bf Corollary 2.2.9.} {\it Let $P=R/S$, $Q=V/W$,
where $R$, $S$, $V$, $W$ are polynomials over $K$, (V, W, monic,
(R,S)$=1$, (V,W)$=1$), $R$, $V$ having degree 2, $\deg W=3$. Let
$P'$ have a zero $c$ of multiplicity $2$, and let $V'-P(c)W'$ have
a zero of multiplicity $2$. If $P(c)\neq Q(d)$, then (\ref{eq1})
has no non constant entire solutions.}\ep
\section{Decompositions of meromorphic functions} Let us
have a look at the functional equation (\ref{eq1}) again, $f$ and
$g$ now being meromorphic functions in all of $K$ (resp. in
$M_u(d(0,r^-))$) and $P$, $Q$ in $K(x)$.\\ What is new and
has to be precisely considered, is that $f$, $g$ might have poles
in $K$. This is the reason for some differences to the preceding
case. We note:
\begin{itemize}
    \item There are not such cases $1$, $2$, $3$, as in Remark 2.2.1 (for the "analytic case"):
    the fact, that $f$, $g$ might have poles yields more "degrees of freedom"
    for the decomposition (\ref{eq1}), i.e.: $\deg R>\deg S\Rightarrow\deg V>\deg W$
    \item If $\deg V=\deg W$, then if $g$ has a pole at $b \in K,\; Q(g)$ has no pole at $b$; this means, in this
    case we cannot get an estimation of $\widetilde
    N(\rho,g)$ by calculating $\widetilde N(\rho, Q(g))$.
    \item Finally, Remark 2.1.1 tells us, that estimation
    (\ref{eq26}) in the proof of Theorem 2.2.4-2.2.5 turns worse, compared with the analytic case.
    \end{itemize}
The aim of this section is to establish statements along the lines
of Theorem 2.2.4 and 2.2.5 for rational decompositions (\ref{eq1})
of two distinct meromorphic functions $f$, $g$. To begin with, we
repeat a statement of \cite{HY8}, presenting a precise asymptotic
formula for the Nevanlinna function of a rational function
composed with a meromorphic one (this is a generalization of
Theorem 0.4 and the analogue formula in Remark 2.2.1):\medskip\\
{\bf Proposition 3.1.1.} {\it Let $f\in {\mathcal M}(K)\setminus
K$ (resp. $f\in {\mathcal M}_u(d(0,r^-)))$), $L\in K(x)$, where
$L=A/B$, and $A$, $B$ in $ K[x]$, $(A,B)=1$ and $\deg A=k,\deg
B=q$. Then we have
\begin{equation}\label{eq101}
T(\rho, L(f))=\max\{k, q\}T(\rho, f)+O(1)
\end{equation}}\medskip\\
For the proof in the case $f\in {\mathcal M}(K)\setminus K$ we
refer to \cite{HY8}, where a little more general statement is
shown. For $f$ being meromorphic in a disk, the proof is
analogue, since the only non elementary facts used are the
Jensen's Formula and the analogue statement Theorem 0.4 for $L\in
K[x]$.\medskip\\Thus we infer the same asymptotic formula for
$T(\rho, g)$ as in the analytic case (Remark 2.2.1):\medskip\\{\bf
Proposition 3.1.2.} {\it If $f\in {\mathcal M}(K)\setminus K$
(resp. $f\in {\mathcal M}_u (d(0,r^-))$), and $f$, $g$ satisfy
(\ref{eq1}), then
\[
q T(\rho, g)=p T(\rho, f)+O(1)
\]}\ep\medskip\\
{\bf Definitions and Notation 3.1.3.} In this section we
distinguish following cases with respect to the degrees of $R$,
$S$, $V$ and $W$ and assign to each of them a certain rational
number $\Lambda(P,Q,f,g)$:
\\{\bf Case 1} $v=w$:            $\;\;\;\quad\quad\Lambda(P,Q,f,g):=\frac{p}{q}$,
\\{\bf Case 2} $v<w$, $r\geq s$: $\Lambda(P,Q,f,g):=\min\{\;\gamma (R),\;\frac{p}{q}\}$,
\\{\bf Case 3} $v>w$, $r\leq s$: $\Lambda(P,Q,f,g):=\min\{\;\gamma (S),\;\frac{p}{q}\}$,
\\{\bf Case 4} $v>w$, $r>s$:     $\Lambda(P,Q,f,g):=\min\{\;\gamma (S)+1,\;\frac{p}{q}\}$,
\\{\bf Case 5} $v<w$, $r<s$:     $\Lambda(P,Q,f,g):=\min\{\;\gamma (R)+1,\;\frac{p}{q}\}$,\\
where for $L \in K[x]$, $\gamma(L)$ denotes the number of distinct
zeros of $L$ in $K$. \\$\Lambda$ arises in following estimation of
$\widetilde N(\rho,g)$ by $T(\rho, f)$:\medskip\\{\bf Proposition
3.1.4.} {\it If $f,\; g\in\mathcal{M}(K)$ (resp. $\mathcal
{M}(d(0,r^-))$) satisfy (\ref{eq1}), then we have
\begin{equation}\label{eq102}
\widetilde N(\rho, g)\leq \Lambda(P,Q,f,g) T(\rho, f)+O(1)
\end{equation}}
The proof is given in section 4.\\ Similarly to Theorem 2.2.4
and 2.2.5 we state now\medskip\\{\bf Theorem 3.1.5.} {\it Let $f$,
$g$  $\in {\mathcal M} (K)\setminus K$, let $P$, $Q$ $\in K(x)$
satisfy condition (M), and let $\Theta(P):=\sum_{i=1}^k (s_i-2)>0
$, i.e., at least one zero $c_j$ of $P-P(c_j)$ has multiplicity
greater than 2. If $f$ and $g$ are solutions to equation
(\ref{eq1}), then, \[ \widetilde N(\rho, g)\geq \left(
\frac{q\Theta(P)-p(k\gamma(W)+1)}{\chi^t q} \right)T(\rho, f)
+\log \rho+O(1)
\]}\medskip\\{\bf Theorem 3.1.6.}
{\it Let $f$, $g$  ${\mathcal M}_u(d(0,r^-))$, let $P$, $Q$ $\in
K(x)$ satisfy condition (M), and let $\Theta(P):=\sum_{i=1}^k
(s_i-2)>0 $, i.e., at least one zero $c_j$ of $P-P(c_j)$ has
multiplicity greater than $2$. If $f$ and $g$ are solutions to
equation (\ref{eq1}), then, \[ \widetilde N(\rho, g)\geq \left(
\frac{q\Theta(P)-p(k\gamma(W)+1)}{\chi^t q} \right)T(\rho, f)
+O(1)
\]}
The proofs can be found in section 4.\medskip\\{\bf Corollary
3.1.7.} {\it Let $f$, $g$  $\in {\mathcal M} (K)\setminus K$, let
$P$, $Q$ $\in K(x)$ satisfy condition (M). If $f$ and $g$ are
solutions to equation (\ref{eq1}), then, we have
\[
q\Theta(P)<p(k\gamma(W)+1)+q\Lambda(P,Q,f,g)
\]}\medskip\\{\bf Corollary 3.1.8.}
{\it Let $f$, $g$  ${\mathcal M}_u(d(0,r^-))$, let $P$, $Q$ $\in
K(x)$ satisfy condition (M). If $f$ and $g$ are solutions to
equation (\ref{eq1}), then, we have
\[
q\Theta(P)\leq p(k\gamma(W)+1)+q\Lambda(P,Q,f,g)
\]}
\medskip\\{\it Proof of the Corollaries 3.1.6-3.1.7.} Both follow
from Theorem 3.1.5 resp. Theorem 3.1.6 by Proposition 3.1.4 (i.e.
the asymptotic formula for $\widetilde N(\rho,g)$) and the growth
of the $\log\rho$- term.\ep
\section{The Proofs}
In this section, we give a unified proof of Theorems 2.2.4, 2.2.5,
3.1.5 and 3.1.6. At first we show that (\ref{eq102}) holds
true:\medskip\\ {\it Proof of Proposition 3.1.4.} In any case we
have $\widetilde N(\rho,g)\leq \frac{p}{q}T(\rho,f)+O(1)$ which
immediately follows from Proposition 3.1.2. In certain cases,
basic considerations improve this asymptotic formula:
\\{\bf Case 1.} As mentioned in the beginning of section 3, in this case poles
    of $g$ are cancelling in $Q(g)$, so no better result for $T(\rho, f)$ than
    the one from above can be achieved.
\\{\bf Case 2.} Taking the reciprocal value of (\ref{eq1}) we
    see that
    \[
    \widetilde N(\rho,g)+\widetilde Z(\rho,V(g))=\widetilde Z(\rho,R(f)),
    \]
    because $w>v$ means that $1/Q(g)$ has a pole if and only if
    $g$ has a pole or $V(g)$ has a zero. Likewise we have $1/P(f)$
    has a pole if and only if $S(f)$ has a zero, since a pole of
    $f$ implied a zero of $1/P(f)$. Using Proposition 3.1.2 we derive
    \[
    \widetilde N(\rho,g)\leq \widetilde Z(\rho,R(f))\leq \gamma (R)T(\rho,f)+O(1)
    \]
\\{\bf Case 3.} Can be worked out like the 2. Case. Indeed, by taking the
    reciprocal value of (\ref{eq1}), $R$ and $S$ merely change their roles.
\\{\bf Case 4.} Obviously any pole of $P(f)$ either is a pole of $f$
    or a zero of $S(f)$, similarly any pole of $Q(g)$ either is a pole
    of $g$ or a zero of $W(g)$, thus we infer
    \[
    \widetilde N(\rho,g)+\widetilde Z(\rho,W(g))=\widetilde
    N(\rho,f)+ \widetilde Z(\rho,S(f)),
    \]
    which means
   \[
   \widetilde N(\rho,g)\leq(1+\gamma(S))T(\rho,f)+O(1)
   \]
\\{\bf Case 5.} Taking the reciprocal value of (\ref{eq1}) we
    conclude as in the preceding case, with the roles of $R$, $S$
    exchanged.
\ep\medskip\\
{\it Proof of Theorems 2.2.4, 2.2.5, 3.1.5 and 3.1.6.} First, let
us suppose $f$ and $g$ to be in ${\mathcal A}(K)\setminus K$
(resp. ${\mathcal A}_u(d(0,r^-))$) and assume $\pi=0$. By Lemma
2.2.2 we have certain decompositions (\ref{eq3}) and (\ref{eq4})
for any fixed $i$, thus by inserting $f$ and $g$ therein we derive
by means of (\ref{eq1})
\begin{equation}\label{eq5}
P(f)-P(c_i)=(f-c_i)^{s_i}R_i(f)=\frac{1}{W(g)}\prod_{j=1}^q
(g-b_{i,j})=Q(g)-P(c_i)
\end{equation}
Applying the second Nevanlinna Theorem N to $g$ we derive
\begin{equation}\label{eq6}
(qk-1)T(\rho, g)\leq \sum_{i=1}^k\sum_{j=1}^q\widetilde
Z(\rho,g-b_{i,j})+\widetilde N(\rho, g)-\log
\rho+O(1)\quad\quad(\rho\rightarrow\infty)
\end{equation}
wherein of course $\widetilde N(\rho, g)=0$; for $g$ we may assume
$g(0)\neq b_{i,j}$ whenever $(i,j)\in\{1,\dots
k\}\times\{1,\dots,q\}$.\\ By means of (\ref{eq4}) we easily
obtain for any fixed $i$
\begin{equation}\label{eq7}
    \widetilde Z(\rho,(Q(g)-P(c_i))W(g))=\sum_{j=1}^q\widetilde
Z(\rho,g-b_{i,j})
\end{equation}
Inserting $g$ in (\ref{eq7}) and in (\ref{eq6}) yields
\begin{equation}\label{eq8}
(qk-1)T(\rho, g)\leq \sum_{i=1}^k\widetilde
Z(\rho,(Q(g)-P(c_i))W(g))-\log \rho+O(1)
\end{equation}
By (\ref{eq2}) we know that $T(\rho,g)=\frac{p}{q}T(\rho,f)+
O(1)$, furthermore by Lemma 2.2.3 we can write
\[
R_i(x)=\frac{A_i(x)}{B_i(x)},\quad \deg(A_i)\leq p-s_i
\]
and obviously $\widetilde N(\rho, B_i(f))=\widetilde N(\rho,
f)=0$, thus by Lemma 2.2.3,
\[
\widetilde Z(\rho, R_i(f))=\widetilde Z(\rho, A_i(f))+\widetilde
N(\rho, B_i(f))\leq (p-s_i)T(\rho,f)+O(1)
\]
For any $i$, any zero of $W(g)$ is a pole of $R_i(f)$ of same
order: this follows from equation (\ref{eq5}) and the fact that
$R$, $S$ have no common zeros. Therefore $\widetilde Z(\rho,
R_i(f)W(g))=\widetilde Z(\rho, R_i(f))$. Thus we can finally
estimate any term of the sum on the right side of (\ref{eq8}) by
\begin{eqnarray}\label{eq26}
&\;&\widetilde Z(\rho,(Q-P(c_i))W(g))=\widetilde Z(\rho,
(f-c_i)^{s_i}R_i(f)W(g) )\leq \\\nonumber &\leq& \widetilde
Z(\rho, (f-c_i))+\widetilde Z(\rho, R_i(f))\leq\\\nonumber&\leq&
T(\rho,f)+(p-s_i)T(\rho,f)+O(1)
\end{eqnarray}
and (\ref{eq40}), (\ref{eq24}) easily follow. \\ It
immediately follows (consider the growth of the $\log \rho$ term)
that for $f\in {\mathcal A}(K)\setminus K$ we have $qk<p$ in
Theorem 2.2.4. Moreover, by Theorem 0.3 we receive $qk\leq p$ in
Theorem 2.2.5.\\ Let now $\pi\neq 0$: By Lemma $\Pi_1$ and
$\Pi_2$ we see that we may apply the Nevanlinna Theorem N in the
same way for $g_t$ having ramification index $0$. So we may write
in the same way as we derived for characteristic $0$:
\[
0\geq (\frac{kq-p}{\chi^tq}) T(\rho, f_t)+\log \rho+O(1),
\]
since the numbers $k$, $p$, $q$ are the same for $P_t, Q_t$. Due
to Lemma $\Pi_2$ the ramification index of $f$ and $g$ are equal
and we immediately get $T(\rho, f_t)=\frac{T(\rho, f)}{\chi^t}$
which finishes the proof of the Theorems 2.2.4, 2.2.5.\\
The proof of the Theorems 3.1.5, 3.1.6
is similar to the one of Theorem 2.2.4 and Theorem 2.2.5:\\
Suppose $\pi=0$. Obviously, formulas (\ref{eq5}), (\ref{eq6}),
(\ref{eq7}) for $f$ and $g$ hold true. Since $g$ might have poles,
instead of (11) we must write now
\begin{equation}\label{eq200}
(qk-1)T(\rho, g)\leq \sum_{i=1}^k\widetilde
Z(\rho,(Q(g)-P(c_i))W(g))+\widetilde N(\rho, g)-\log \rho+O(1)
\end{equation}
Now we need to estimate (\ref{eq7}): For any $i=1,\dots,k$, we
receive
\[
\widetilde Z(\rho, R_i(f))\leq \widetilde Z(\rho,
A_i(f))+\widetilde N(\rho, f)\leq (p-s_i+1)T(\rho,f)+O(1),
\]
thus,
\begin{eqnarray}\label{eq103}
&\;&\widetilde Z(\rho,(Q(g)-P(c_i))W(g))=\widetilde Z(\rho,
(f-c_i)^{s_i}R_i(f)W(g) )\leq \\\nonumber &\leq& \widetilde
Z(\rho, (f-c_i))+\widetilde Z(\rho, R_i(f))+\widetilde Z(\rho,
W(g))\leq\\\nonumber&\leq& T(\rho,f)+(p-s_i+1)T(\rho,f)+\gamma
(W)T(\rho,g)+O(1)=\\\nonumber&=&\left(p-(s_i-2)+\frac{p
\gamma(W)}{q}\right)T(\rho,f)+O(1)
\end{eqnarray}
Summing up over $i$ we derive
\[
\sum_{i=1}^k \widetilde Z(\rho,(Q(g)-P(c_i))W(g))\leq
(kp-\Theta(P)+\frac{kp\gamma(W)}{q})T(\rho,f)+O(1)
\]
And this estimation put into (\ref{eq200}) yields the maintained
result.\\ If $K$ has characteristic $\pi\neq 0$, the proof
is similar to the one of the Theorems 2.2.4-2.2.5 in the same
situation. Note that $\widetilde N(\rho,g)=\widetilde
N(\rho,g_t)$.\ep \medskip\\{\bf Acknowledgement}\medskip\\ I want
to thank Professor Alain Escassut for suggesting this topic and
for giving a lot of useful hints.


\begin{thebibliography}{1}

\bibitem{B1}
{\sc A.~Boutabaa}, {\em On some {$p$}-adic functional equations}, in $p$-adic
  functional analysis (Nijmegen, 1996), vol.~192 of Lecture Notes in Pure and
  Appl. Math., Dekker, New York, 1997, pp.~49--59.

\bibitem{BE2}
{\sc A.~Boutabaa and A.~Escassut}, {\em Applications of the {$p$}-adic
  {N}evanlinna theory to functional equations}, Ann. Inst. Fourier (Grenoble),
  50 (2000), pp.~751--766.

\bibitem{BE3}
{\sc A.~Boutabaa and A.~Escassut}, {\em Urs and ursims for {$p$}-adic
  meromorphic functions inside a disc}, Proc. Edinb. Math. Soc. (2), 44 (2001),
  pp.~485--504.

\bibitem{BE4}
\leavevmode\vrule height 2pt depth -1.6pt width 23pt, {\em Nevanlinna theory in
  characteristic {$p$} and applications}, in Analysis and applications---ISAAC
  2001 (Berlin), vol.~10 of Int. Soc. Anal. Appl. Comput., Kluwer Acad. Publ.,
  Dordrecht, 2003, pp.~97--107.

\bibitem{E5}
{\sc A.~Escassut}, {\em Analytic elements in {$p$}-adic analysis}, World
  Scientific Publishing Co. Inc., River Edge, NJ, 1995.

\bibitem{EY6}
{\sc A.~Escassut and C.-C. Yang}, {\em The functional equation {$P(f)=Q(g)$} in
  a {$p$}-adic field}, J. Number Theory, 105 (2004), pp.~344--360.

\bibitem{G7}
{\sc F.~Gross}, {\em On the equation $f^n+g^n=1$}, Bull. Amer. Monthly, 73
  (1966), pp.~1093--1096.

\bibitem{HY8}
{\sc P.-C. Hu and C.-C. Yang}, {\em Meromorphic functions over
  non-{A}rchimedean fields}, vol.~522 of Mathematics and its Applications,
  Kluwer Academic Publishers, Dordrecht, 2000.

\bibitem{T9}
{\sc N.~Toda}, {\em On the functional equation {$\sum
  \sp{p}\sb{i=0}\,a\sb{i}f\sb{i}\sp{n\sb{}i}=1$}}, T\^ohoku Math. J. (2), 23
  (1971), pp.~289--299.

\end{thebibliography}
\end{document}